\documentclass[11pt]{article}
\usepackage{amsmath}
\usepackage{amssymb}
\usepackage{amsfonts} 
\usepackage{amstext} 
\usepackage{amsthm} 
\usepackage{euscript} 
\usepackage{eufrak}

\newcommand{\Z}{\ensuremath{\mathbb{Z}}} 
\newcommand{\N}{\ensuremath{\mathbb{N}}}
\newcommand{\Q}{\ensuremath{\mathbb{Q}}}
\newcommand{\R}{\ensuremath{\mathbb{R}}} 
\newcommand{\C}{\ensuremath{\mathbb{C}}}

\newcommand{\cM}{\ensuremath{\mathcal{M}}} 
 
\newcommand{\cP}{\ensuremath{\mathcal{P}}}

\newcommand{\cA}{\ensuremath{\mathcal{A}}}

\begin{document} 
\sloppy 
\begin{center} 
{\Large \bf Vanishing of the Contact Homology of  \\  
Overtwisted Contact 3-Manifolds} \\ 

\vspace{.3in} 
Mei-Lin Yau  \\  
(with an appendix by Yakov Eliashberg) \\ 
   \vspace{.2in} 
Department of Mathematics \\  Michigan State University \\ 
East Lansing, MI 48824 \\ 
Email: yau@math.msu.edu
\end{center} 
\vspace{.2in} 
\begin{abstract} 
We give a proof of, for the case of contact structures defined by 
global contact 1-forms, a theorem stated by 
Eliashberg that for 
any overtwisted contact structure 
on a closed 3-manifold, its contact homology is 0. A different proof 
is also outlined in the appendix by Yakov Eliashberg. 
\end{abstract}


\newtheorem{cond}{Condition}[section]
\newtheorem{defn}{Definition}[section]
\newtheorem{exam}{Example}[section]
\newtheorem{lem}{Lemma}[section]
\newtheorem{cor}{Corollary}[section]
\newtheorem{main}{Main Theorem}
\newtheorem{theo}{Theorem}[section]
\newtheorem{prop}{Proposition}[section]
\newtheorem{rem}{Remark}[section]
\newtheorem{notn}{Notation}[section] 
\newtheorem{fact}{Fact}[section] 





\section{Introduction} \label{intro} 

A contact structure $\xi$ on an odd dimensional manifold is a 
nowhere integrable hyperplane distribution. If $\xi$ is coorientable 
then it is defined by a global 1-form $\alpha$, i.e., $\xi:=\ker\alpha$. 
$\alpha$ is called a defining {\em contact 1-form} of $\xi$. 
{\em In this article all contact structures which we consider are 
defined by global  contact 1-forms, 
and all manifolds are closed and 
orientable.} A contact manifold $(M,\xi)$ consists of a 
$(2n-1)$-dimensional
manifold $M$ and a contact structure $\xi$ on $M$.  Write
$\xi:=\ker\alpha$ where $\alpha$ is a contact 1-form  defining
$\xi$, then $(M,\xi)$ has a natural orientation  defined by the
volume form $\alpha\wedge (d\alpha)^{n-1}$.

A contact structure $\xi$ on a 3-manifold $M$ is 
{\em overtwisted} 
if there exists an embedded disc $D\subset M$ such that 
$T(\partial D)\subset \xi|_{\partial D}$, 
$T_zD\not\subset\xi_z$ for all $z\in \partial D$, 
and $\xi\cap TD$ defines a singular foliation with exactly one 
singular point (which has to be elliptic). 
It is proved by 
Eliashberg \cite{E1} that every homotopy class of plane
distributions  of a 3-manifold has a unique (up to a contact isotopy) 
overtwisted contact representative. On the other hand, many
3-manifolds possess non-overtwisted contact structures, 
making the classification of contact manifolds a very subtle 
question.

Based on the introduction of pseudoholomorphic curves 
in symplectic manifolds by Gromov \cite{Gr}, Eliashberg and Hofer 
in the mid 90's introduced 
Contact Homology Theory (\cite{E2} 
\cite{EGH}, see also \cite{B} for the Morse-Bott version)  
to provide Floer-type invariants $H\Theta(\xi)$, the {\em contact 
homology}, for contact manifolds.  
Here let us describe very briefly the construction of contact 
homology, which is in fact an algebra. Readers are suggested to
consult  
\cite{E2}, \cite{EGH} and \cite{B} 
for more detail. First of all, each contact 1-form 
$\alpha$ associates a unique vector field $R=R_\alpha$ transversal 
to $\xi:=\ker\alpha$, $R$ is called the {\em Reeb} vector 
field associated to $\alpha$ and is defined by 
\[ 
d\alpha (R,\cdot )=0, \quad \alpha (R)=1. 
\] 
The {\em symplectization} $(\R\times M,d(e^t\alpha))$ ($t\in \R$) 
is equipped with an $\alpha$-{\em admissible} 
(see Section \ref{ch-alg2}) 
almost complex complex structure. 
The contact homology is a homology whose complex is generated 
by {\em good} (see Section \ref{ch-alg}) 
periodic Reeb trajectories, and whose boundary 
operator $\partial$ is defined by counting in  
$\R\times M$ 
one-dimensional moduli of pseudoholomorphic 
curves of {\em genus 0} with finite $d\alpha$-energy 
which converge asymptotically to a single Reeb orbit as $t\to\infty$ 
and an arbitrary number of Reeb orbits as $t\to -\infty$. 
The resulting homology depends only on the isotopy class of 
contact structures (\cite{E2},\cite{EGH}).

The main purpose of this paper 
is to give a proof to the following theorem 
stated by Y. Eliashberg 
\cite{E2} in the case where the overtwisted contact structure is
defined by a  global contact 1-form.

\begin{theo}[Eliashberg, \cite{E2}] \label{ot}
If $\xi$ is an overtwisted contact structure on a 3-manifold $M$, 
then $H\Theta(\xi)=0$. 
\end{theo}

Note that in 
Hofer's proof of Weinstein Conjecture \cite{H} 
for overtwisted contact 3-manifolds he showed there 
exists a contractible Reeb orbit in 
every overtwisted contact 3-manifold, and 
such a Reeb orbit must be the asymptotic boundary of a 
pseudoholomorphic plane with finite energy. {\em If}  
this pseudoholomorphic plane is the {\em only}  
pseudoholomorphic curve that converges to the said Reeb orbit 
as $t\to\infty$, then the contact homology of the overtwisted 
contact structure is 0. 

In this paper our proof of Theorem \ref{ot} is based on 
the classification of overtwisted contact structures by 
Eliashberg \cite{E1},  
open book representations of contact manifolds 
(see Theorem \ref{TWG}) 
by Thurston and Winkelnkemper \cite{TW} and 
Giroux \cite{G1},  
a construction of an overtwisted contact structure from a 
trivial Dehn surgery inspired by Geiges \cite{Ge}, as well 
as some techniques and conclusions of some $S^1$-invariant 
moduli of pseudoholomorphic curves from \cite{Y2} and 
finally, the fact that contact homology is independent of 
the choices of a contact 1-form and an admissible almost 
complex structure in the construction of $H\Theta(\xi)$.

Here is an outline of the paper: Section \ref{obook} 
consists of some background on open books. In Section 
\ref{ot-dehn} we first construct a contact 1-form $\alpha$ 
from an open book then "twist" it along a contractible 
Reeb orbit by a trivial Dehn surgery to get a new contact 
1-form $\alpha'$. It is then proved that $\xi':=\ker\alpha'$ 
is overtwisted and yet is homotopic to $\xi$ as plane 
distributions. In particular we get a special contractible 
Reeb orbit $t_x$ which will be proved in Section \ref{ch} to 
satisfy the equation $\partial t_x=\pm 1$. 
Section \ref{ch-alg} consists of brief definitions of 
contact complex and contact homology, as well 
as some discussion on cylindrical contact homology. 
In Section \ref{ch} we first show that energy and homotopy 
constraints severely limit the types of pseudoholomorphic curves 
asymptote $t_x$ at positive infinity. Such holomorphic curves 
must be finite energy planes. Then by using methods 
similar to \cite{Y2} we prove that, modulo free $\R$-actions, 
the algebraic number of 
such holomorphic planes is equal to the algebraic number of 
certain gradient trajectories of a Morse function, hence is 
$\pm 1$. Thus $\partial t_x=\pm 1$ and the contact homology of 
$\xi'$ is 0. 

\vspace{.3in} 
\noindent 
{\bf Acknowledgments:} 

After this paper is written we were kindly reminded \cite{E3}  
that a proof to Theorem \ref{ot} was already known to Eliashberg 
before his talk at ICM-Berlin 1998 \cite{E2}, albeit not written. 
Eliashberg also generously offered us an outline of his proof 
\cite{E3} which we include in the appendix of this paper.  


\section{Open books and contact structures} \label{obook}

For any surface $S$ and any $\psi\in\text{Diff}(S)$ 
we denote by $S_{\psi}$ the mapping torus 
\begin{gather*}  
S\times [0,1]/\sim, \\  
(\psi (x),0)\sim (x,1). 
\end{gather*}

An {\em open book decomposition} of a connected closed orientable 
3-manifold $M$ consists of a 1-dimensional 
submanifold $B$ (a link in $M$), called the {\em binding}, 
and a fibration 
$\pi :M\setminus B \to S^1$ with fibers connected embedded surfaces 
with boundary $B$. The fibers are called {\em pages}.  
 In this paper we  
assume that there is a tubular neighborhood $B\times D^2$ of 
$B$ so that $\pi$ restricts to the normal angular coordinate of 
$B=B\times \{ 0\}$ in $B\times D^2$. Then 
\[  
M=\Sigma _\phi \cup_{id}(B\times D^2),  
\]  
where $\Sigma$ is an orientable surface 
with boundary $\partial \Sigma\cong B$ and 
$\phi\in\text{Diff}^+(\Sigma, \partial \Sigma)$ an 
orientation-preserving diffeomorphism with $\phi=id$ near
$\partial \Sigma$. Note that 
$\phi$ is unique up to isotopy. 
The pages of $\pi :M\setminus B\to S^1$ are diffeomorphic to
$\Sigma$.

We will need the following important result by 
Thurston and Winkelnkemper \cite{TW} and Giroux \cite{GM} 
concerning contact structures on 3-manifolds. 

\begin{theo} \label{TWG} 
Each
open book associates a  unique up to isotopy contact structure and
conversely,  each contact structure is supported by an open book
unique  up to {\em positive stabilizations}. 
\end{theo} 

We will not go into the detail here but 
point out that by applying positive stabilizations several times
if necessary we may assume that $\partial\Sigma\cong B$ is 
{\em connected}. 

With Theorem \ref{TWG} in hand we can start with any open book 
$(\Sigma, \phi)$ and alter the corresponding contact structure $\xi$ 
by a trivial Dehn surgery to get a new contact structure $\xi'$ 
homotopic to $\xi$ as plane distribution, as we will do in the 
next section.

\section{Overtwisted contact form from trivial Dehn surgery} 
\label{ot-dehn}

In this section we will construct an overtwisted contact 
structure $\xi':=\ker\alpha'$ 
that is homotopic to a given contact structure 
$\xi:=\ker\alpha$ as plane distributions. We sketch the 
idea of the 
construction of $\alpha'$ here before going into the detail. 

Start with a 
contact 1-form $\alpha$ (following \cite{TW}) 
associated to a given open 
book $(\Sigma,\phi)$ such that 
$\alpha$ has a pair of contractible 
Reeb orbits $t_x$ and $t_y$ associated to a pair of 
birth-death type of critical points $x,y$ of a smooth function 
$K$ on $\Sigma$, such that $x$ is a saddle point and $y$ is a
local  minimum point. Following the construction in \cite{Ge}  
of contact 1-forms under Dehn surgeries, 
we apply a trivial Dehn surgery to $t_y$. 
In particular, we cut out a tubular neighborhood of $t_y$  
and glue it back, identifying the two boundaries  
by using the gluing matrix $\begin{bmatrix}-1&0\\0&-1\end{bmatrix}$. 
The resulting 3-manifold is the 
same (up to an orientation-preserving diffeomorphism) while  
the contact 1-form $\alpha$, after being modified near the 
boundaries, glued back to a new contact 1-form $\alpha'$. 
The new contact structure $\xi'$ is then shown to be overtwisted 
and yet is homotopic to $\xi$ as plane distributions. 
Now start the construction.

Given an open book $(\Sigma,\phi)$ with 
$\partial\Sigma \cong B$ 
connected, we can define an associated  
contact 1-form $\alpha$ on 
$M=\Sigma_\phi\cup_{id}B\times D^2$ as follows.

Let $F\subset\Sigma$ be a collar of $\partial\Sigma$ 
such that $\phi |_F=id$. Let $(q,p)$ be coordinates of 
$F=[q_-,q_+]\times S_p^1$ so that $\partial\Sigma$ is identified 
with $\{q_+\}\times S_p^1$. We may assume that 
$0\ll q_-\ll q_+$. 
Let $(\rho,t)$ be the polar coordinates of 
the $D^2$-factor of $B\times D^2\cong S^1\times D^2$ 
so that $D^2=\{ \rho\leq 1\}$. Let 
$p$ be the coordinate of $S^1$. Note that $\rho$ can be viewed 
as a smooth function of the coordinate $q$ of $F$, and 
$\frac{d\rho}{dq}<0$. Also, 
$S^1_t$ acts on $F_\phi\cup (B\times D^2)$ via 
rotations in the $t$-direction. It fixes $B$ and acts 
freely on the complement of $B$ in $F_\phi\cup B\times D^2$. 
Let $V$ denote the orbit space and let 
\begin{equation}  \label{pit} 
\pi_t:F_\phi\cup (B\times D^2)\to V
\end{equation}  
denote the corresponding projection.

Let $K\geq 0$ be a smooth function on $M$ 
such that 

\begin{enumerate} 
\item $K|_{(\Sigma\setminus F)_\phi}\gg 1$ is a constant; 

\item $K$ is $S^1_t$-invariant on $F_\phi\cup (B\times D^2)$; 
$\bar{K}:=(\pi_t)_*K$  is 
a Morse function on $int(V)$ the interior of $V$ with 
$|dK|\sim 0$ on $F$, $K=1$ on $\partial\Sigma$; 

\item $\bar{K}$ has precisely two critical 
points on $int(V)$: one 
saddle point $x$ and one local minimum point $y$; 
$x$ and $y$ are of birth-death type, i.e., there is 
only one gradient trajectory of $\bar{K}$ 
that connects $x$ and $y$; 

\item $\bar{K}(x)>\bar{K}(y)>1$ 
and $\bar{K}(x)\sim \bar{K}(y)\sim 1$;  

\item in the interior of $F$ choose a 
small disc neighborhood 
\[ D_y=\{ (r,\theta)\mid r\leq \delta\} ,\quad \delta>0 \text{ a
constant},  
\]  
with center $y$,  here $(r,\theta)$ are 
the polar coordinates of $D_y$, then 
$\bar{K}|_{D_y}$ depends only on $r$; 

\item $K|_{B\times D^2}$ depends only on $\rho$, 
$\frac{dK}{d\rho}>0$ on $0<\rho\leq 1$, $K(\rho)=h(\rho)\rho^2$ 
for some smooth function $h$ depending only on $\rho$. 
\end{enumerate}

Let $\beta$ be a 1-form on $\Sigma$ such that 
$d\beta$ is a symplectic 2-form on $\Sigma$ with 
\[ 
\beta=qdp \text{ near }\partial\Sigma, \qquad 
\beta=r^2d\theta \text{ on } D_y . 
\]

On $B\times D^2$ we also consider 
a smooth function $Q>0$ depending only on $\rho$ such that 
\begin{gather*}  
Q=q \text{ near }\rho=1,  \qquad \frac{dQ}{d\rho}(0)=0, \\ 
\frac{dQ}{d\rho}>0\text{ and } 
Q\frac{dK}{d\rho}-K\frac{dQ}{d\rho}>0 \text{ for }\rho>0. 
\end{gather*} 
 
Then 
\begin{equation}  \label{alpha} 
\alpha :=\begin{cases} 
(1-t)\beta +t\phi^*\beta +Kdt  \quad 
\text{ on }\Sigma_\phi  \\ 
Q(\rho)dp+K(\rho)dt \hspace{.5in} \text{ on } B\times D^2 
\end{cases} 
\end{equation}  
is a contact 1-form on $M$ provided that $K|_{\Sigma\setminus F}$ 
is a large enough constant. Denote $\xi:=\ker\alpha$.

We have two contractible simple Reeb orbits of $\alpha$: 
\[ 
t_x:=\{ x\}\times S^1_t ,\quad t_y:=\{ y\}\times S^1_t. 
\] 
Both $t_x$ and $t_y$ are oriented by the vector field 
$\partial_t$.

Now define  
a contact 1-form $\alpha'$ on $M$:  
\begin{equation}  
\alpha'=\begin{cases}   \alpha \hspace{1.2in} \text{ on
}M\setminus  (D_y)_\phi, \\ h_1(r)dt+h_2(r)d\theta \quad \text{ on }
(D_y)_\phi. 
\end{cases} 
\end{equation}  
where $h_1$, $h_2$ are smooth functions of $r$ satisfying 
\begin{enumerate} 
\item $h_1(r)=-1$ and $h_2(r)=-r^2$ for $r<\epsilon$; 
\item $h_1(r)=K(r)$ and $h_2(r)=r^2$ for
$\epsilon\leq r\leq\delta-\epsilon$; 
\item $h'_1(r)h_2(r)-h'_2(r)h_1(r)>0$ for $0<r\leq\delta
-\epsilon$;  
\end{enumerate} 
where $0<\epsilon\ll \frac{\delta}{2}$ is a constant. 
The third condition above is to ensure that $\alpha '$ is a contact 
form on $(D_y)_\phi$. Note that $\alpha'$ has two special Reeb 
orbits 
\[ 
t_x=\{ x\}\times S^1_t,\quad \bar{t}_y=\{y\}\times S^1_t 
\] 
The notation $\bar{t}_y$ represents the curve $t_y$ but with the 
reversed orientation, i.e., the orientation given by $-\partial_t$. 
Let $\xi':=\ker \alpha'$. 

\begin{lem} 
$\xi'$ is an overtwisted contact structure. 
\end{lem} 

\begin{proof} 
Let $\ell$ be the gradient trajectory of $-\bar{K}$ that goes from 
$x$ to $y$. Then $\ell\times S^1_t$ is a homotopy between 
$t_x$ and $t_y$. We have 
\[ 
\int _{t_y}\alpha '<0<\int_{t_x}\alpha ' 
\] 
so there exists a point $z\in \ell$ such that 
$\int_{t_z}\alpha'=0$. Since $\alpha'$ is $S^1_t$-invariant on 
$\ell\times S^1_t$ we conclude that 
\[ t_z \text{ is a Legendrian curve of }\xi'. 
\] 
Since $t_z$ is homotopic to $t_x$, $t_z$ is 
contractible. Moreover, $t_z$ is contained in a tubular
neighborhood of the binding $B$ and its winding number with 
$B$ is $\pm 1$, so $t_z$ is spanned by an embedded disc in 
a tubular neighborhood of $B$. We can find an overtwisted 
disc spanning $t_z$ as follows. 

Note that $h_1$ vanishes at $z$, and that 
$\partial_t \pitchfork \xi '$ on $F_\phi\cup (B\times D^2)$ 
except where $h_1=0$. Recall that $V$ is the projection 
of $F_\phi\cup (B\times D^2)$ via $\pi_t$. 
Let $\gamma\subset V\setminus \{ h_1(r)<0\}$ be an   
embedded smooth path 
such that 
\begin{enumerate} 
\item $z$ is an endpoint of $\gamma$, 
\item $\gamma|_{D_y}$ is transversal 
to $\partial_r$, and 
\item 
$\dot{\gamma}\parallel \partial_\rho$ near $B$. 
\end{enumerate} 
Define 
$D_z:=\pi_t^{-1}\gamma$. Then $D_z$ is an embedded 
smooth spanning disc of
$t_z$. Moreover, since $\xi'|_{t_z}=\text{Span}(\partial_t,\partial_r)$, 
$D_z\pitchfork \xi '$ along $t_z$ by the 
second condition above, $\partial_t$ is tangent to 
$D_z\setminus\{ \rho =0\}$ hence $\xi'\cap TD_z$ has only 
one singular point and the singular point is elliptic. 
So $D_z$ is an overtwisted disc.  
\end{proof} 

\begin{lem} 
$\xi'$ and $\xi$ are homotopic as plane distributions. 
\end{lem} 
 
\begin{proof} 
It is enough to show that 
 $\alpha$ and $\alpha'$ are homotopic 
as nowhere vanishing 1-forms. Since $\alpha=\alpha'$ on 
$M\setminus (D_y)_\phi$ it is enough to consider a homotopy 
supported in $(D_y)_\phi$. 

For $s\in [0,1]$ define 
\[ 
\alpha_s:=s(1-s)\chi (r)dr+(1-s)\alpha +s\alpha ' \quad 
\text{ on }(D_y)_\phi
\] 
where $\chi(r)\geq 0$ is a smooth function on $r$ such that 
$\chi(0)=0=\chi(\delta)$ and $\chi >0$ away from $r=0,\delta$. 
We have 
\[ 
\alpha_s=s(1-s)\chi(r)dr+\big( (1-s)K(r)+sh_1(r)\big) dt+
\big( (1-s)r^2+sh_2(r)\big) d\theta. 
\] 
It is clear that $\alpha_s$ is nowhere vanishing when $s$ is 
close to $0$ or $1$. Also for every $s$, $\alpha_s$ is
nonvanishing on the region where $\chi$ is positive and 
the region where $r$ is close to 1. 
For $s\neq 0,1$ we have 
\begin{gather*} 
(1-s)K(r)+sh_1(r)=0 \Leftrightarrow h_1(r)=\frac{(s-1)K(r)}{s} \\ 
(1-s)r^2+sh_2(r)=0 \Leftrightarrow h_2(r)=\frac{(s-1)r^2}{s} 
\end{gather*} 
The function $g(s):=\frac{s-1}{s}$, $s\in (0,1)$, is an 
increasing function of $s$ such that $g(s)\to -\infty$ as $s\to
0^+$, and $g(s)\to 0^-$ as $s\to 1^-$. 

Recall that near $r=0$, 
\[ 
\frac{h_1(r)}{K(r)}=\frac{-1}{K(r)}, \quad 
\frac{h_2(r)}{r^2}=\frac{-1}{1} 
\] 
Since $K(r)\neq 1$ for $r$ near $0$,  
for $r$ small enough 
$(1-s)K(r)+sh_1(r)$ and $(1-s)r^2+sh_2(r)$ will not vanish 
simultaneously for any $s\in [0,1]$. Therefore $\alpha_s$ 
is nowhere vanishing for $s\in [0,1]$. Since
$\alpha_0=\alpha$ and $\alpha_1=\alpha '$, $\alpha$ and 
$\alpha'$ are homotopic as nowhere vanishing 1-forms on $M$. 
So $\xi$ and $\xi'$ are homotopic as plane distributions.  

\end{proof}

\section{An outline of contact homology} 
\label{ch-alg} 

In this section we give a brief account on definitions  
of contact complex, contact homology and cylindrical 
contact homology. Readers are referred to 
\cite{E2}\cite{EGH}\cite{B} for more detail. 

\subsection{Contact complex algebra}

Let $(M,\xi)$ be a $(2n-1)$-dimensional closed contact manifold 
with $\xi$ defined by a global contact 1-form $\alpha$. For a 
generic choice of $\alpha$, there are only countably many 
periodic trajectories (including all positive multiple ones) 
of the Reeb vector field $R_{\alpha}$; and these Reeb orbits 
are {\em nondegenerate}, meaning that 
$1$ is not an eigenvalue of their Poincar'{e} return 
map. We call such contact 1-forms {\em regular}.  

\begin{defn}  \label{good} 
A Reeb orbit is said to be {\em bad} 
(see Section 1.2 of \cite{EGH}) if it is an even multiple 
of another Reeb orbit whose Poincar\'{e} return map has the 
property that the total multiplicity of its eigenvalues 
from the interval $(-1,0)$ is odd. A Reeb orbit is {\em good} 
if it is not bad. 
\end{defn} 
We denote by $\cP_{\alpha}$ the set of all good Reeb orbits 
of $\alpha$. Note that $\cP_{\alpha}$ includes all positive 
multiple ones as individual elements.

%
%
%
%

\begin{defn} 
Let $\alpha$ be a regular contact 1-form defining the contact 
structure $\xi$. The contact complex $\Theta (\alpha)$ is then 
defined to be 
the free commutative algebra over $\Q$ (or $\R$, $\C$) 
generated by all elements of $\cP_\alpha$. 
\end{defn} 

\begin{rem} {\rm In \cite{EGH} $\Theta(\alpha)$ is defined 
with coefficients in the algebra $\C[H_2(M)][[t]]$. Here 
we use $\Q$-coefficients for the sake of simplicity. 
}
\end{rem} 

%
%
%
%

\subsection{Contact homology} \label{ch-alg2} 

An almost complex structure $J$ 
on the symplectization 
$(\mathbb{R}\times M,d(e^t\alpha))$  
of $(M,\xi =\ker \alpha )$ 
is said to be $\alpha$-{\em admissible} if 
$J(\partial_t)=R_\alpha$ and 
$J|_\xi:\xi\to
\xi$ on $\xi$ is {\em $d\alpha$-compatible}, i.e., 
$d\alpha (v,Jv) > 0$ for all nonzero $v\in \xi$ and 
$d\alpha (Jv_1,Jv_2)=d\alpha (v_1,v_2)$ for $v_1,v_2\in \xi$.
Note that compatibility property does not depend on the choice
of the defining contact 1-form for $\xi$. 

Let $\gamma$ be a good Reeb orbit. 
We use the following notations: 
\begin{enumerate} 
\item $\Upsilon$:= a 
finite collection of (not necessarily distinct) good Reeb orbits 
of $\alpha$. $\Upsilon$ can be empty.
\item  $|\Upsilon|$:=
the cardinality of $\Upsilon$. 
\item $\tilde{\cM}(\Upsilon,\gamma)$:= the moduli space of 
finite $d\alpha$-energy 
pseudoholomorphic maps from a $(1+|\Upsilon|)$-punctured sphere 
into $\R\times M$ with one puncture goes to $\gamma$ at $t=\infty$ 
and other punctures go to $\Upsilon$ at $t=-\infty$ (see \cite{EGH}). 
\item $\cM (\Upsilon,\gamma)$:=
the union of 1-dimensional components of
$\tilde{\cM}(\Upsilon,\gamma)$.
\end{enumerate}  
Note that $\Upsilon$ was treated 
in \cite{EGH} as an {\em ordered} set, yet {\em here we consider 
$\Upsilon$ an unordered set}. 

Secondly, an admissible almost complex structure $J$ on $\R
\times M$ is $\R$-invariant, hence there is a free $\R$-action 
on $\tilde{\cM}(\Upsilon,\gamma)$.  For generic $\alpha$-admissible 
$J$, 
$\cM(\Upsilon,\gamma)/\R$ consists of finitely many points. 
Note that because of energy constraint (\ref{energy}) there 
are only finitely many choices for $\Upsilon$ such that 
$\tilde{\cM}(\Upsilon,\gamma)\neq\emptyset$. 
Let $\kappa_{\gamma}$ denote the multiplicity of $\gamma$. 
For $C\in \cM(\Upsilon, \gamma)/\R$ we also denote by 
$\kappa_C$ the multiplicity of $C$.

The boundary operator $\partial$ 
of the contact complex $(\Theta(\alpha),\partial)$, 
when applied to $\gamma$, is defined by 
(see \cite{B2}\cite{E2}\cite{EGH} but for a different coefficient
ring) 

\begin{align}   
\partial \gamma & :=\sum_{i=0}^\infty \partial_i\gamma , 
\quad \quad \text{ where } \\ 
\partial_i \gamma & :=\kappa_{\gamma}
\sum_{|\Upsilon|=i}\Big( 
\underset{\dim
\cM(\Upsilon,\gamma)=1}{\sum_{C\in\cM(\Upsilon,\gamma)/\R}}  
\frac{\pm 1}{\kappa_C}\Big) \Upsilon.  \label{coeff} 
\end{align}   

Here $\Upsilon$ denotes the monomial 
$\gamma_1\gamma_2\cdots\gamma_{\Upsilon}$, and the $\pm$ sign 
in (\ref{coeff}) depends on the orientation of 
$C\in\cM(\Upsilon,\gamma)/\R$.

Note that because the action $\cA (\sigma):=\int_{\sigma}\alpha$ 
of any Reeb orbit $\sigma$ of $\alpha$ is bounded from 
below by a positive number independent of $\sigma$, 
$\partial _i\gamma=0$ for all 
$i$ sufficiently large, and the right hand side  of 
(\ref{coeff}) consists of only 
finitely many nonvanishing terms. 
Then extend $\partial$ over 
$\Theta(\alpha)$ according to the Leibnitz rule \cite{EGH}. 
$(\Theta(\alpha), \partial)$ is now a differential algebra.

\begin{theo}[\cite{EGH}] 
$\partial ^2=0$ for regular contact 1-form $\alpha$ and 
a generic $\alpha$-admissible almost complex structure $J$. 
\end{theo} 

\begin{defn} 
The contact homology algebra $H\Theta(\alpha,J)$ of 
the pair $(\alpha,J)$ with $J$ $\alpha$-admissible 
is defined to be the quasi-isomorphism class of 
the differential algebra $(\Theta(\alpha),\partial)$, i.e.,  
$H\Theta(\alpha,J):=\dfrac{\ker\partial}{\text{im}\partial}$.  

\end{defn} 

$H\Theta (M,\xi):=H\Theta(\alpha,J)$ 
is independent of the choices of $(\alpha,J)$ hence  
is an invariant of the contact manifold $(M,\xi)$.

\subsection{Cylindrical contact homology}

Recall the boundary operator 
$\partial=\sum_{i=0}^{\infty}\partial_i$.  
Since $\partial^2=0$ we have 
\[ 
\partial_1^2+\partial_0\partial_2+\partial_2\partial_0=0 
\] 
{\em If} $\partial_1^2=0$ then one can forget the algebraic 
structure of $\Theta(\alpha)$, 
simply thinking it as a free module over $\Q$ generated by 
all good Reeb orbits, and then use $\partial_1:\Theta\to\Theta$ 
as the boundary operator to 
define the {\em cylindrical contact homology} $HC(\xi):=
\dfrac{\ker\partial_1}{\text{im}\partial_1}$ of $(M,\xi)$ as a 
vector space. A sufficient (but not necessary) 
condition for $\partial_1^2=0$ to 
be true is $\partial_0=0$. Note that since $\partial_0$ is defined 
via counting finite $d\alpha$-energy 
pseudoholomorphic planes bounding contractible Reeb orbits at 
$t=\infty$, $\partial_0=0$ holds trivially when there exist no 
contractible Reeb orbits.

\section{Contact homology of overtwisted $\xi'$} \label{ch} 
 
Back to the overtwisted contact 3-manifold $(M,\xi')$ 
and the contact 1-form $\alpha'$ that we constructed 
in Section \ref{ot-dehn}. From Definition \ref{good} it 
is clear that $t_x$ is {\em good} hence is a generator of 
the contact algebra $\Theta(\alpha')$ 
for the contact homology of $\xi'$.

From now on $\Upsilon$ denotes a 
finite collection of (not necessarily distinct) good Reeb orbits 
of $\alpha'$. Note that $\Upsilon$ can be empty. The notations 
$\tilde{\cM}(\Upsilon,t_x)$ and $\cM (\Upsilon, t_x)$ are as 
defined in Section \ref{ch-alg} (with respect to the overtwisted 
contact 1-form $\alpha'$). 
Here we point out a few things about $\tilde{\cM}(\Upsilon,t_x)$ 
which we will need later.

Given $\tilde{u}\in \tilde{\cM}(\Upsilon,t_x)$ we write 
$\tilde{u}=(a,u)$ according to the splitting $\R\times M$  
and let $C\subset M$ denote the image of $u$. {\em Assume that 
$\Upsilon\neq t_x$}, then 
\begin{gather}  \label{pos} 
d\alpha'>0 \text{ on $C$ except at finitely many points of }  C. 
\end{gather} 
Moreover the $d\alpha'$-energy of $\tilde{u}$ is 
\begin{equation} \label{energy} 
\int_C d\alpha'=\int_{t_x}\alpha'-\int_\Upsilon\alpha' >0,  
\end{equation}  
unless when $\Upsilon =t_x$. In the exceptional case, 
$\tilde{\cM}(t_x,t_x)=\text{pt}$, 
the corresponding pseudoholomorphic curve 
is the trivial cylinder $\R\times t_x$. 
Since we are only interested in 
holomorphic curves with positive finite $d\alpha'$-energy, we 
may assume that $\Upsilon$ does not contain $t_x$ and its 
positive iterates.

In the following we will study 
$\tilde{\cM}(\Upsilon, t_x)$ and $\cM(\Upsilon,t_x)$. 
First a few more notations. 

Recall that $F\subset \Sigma$ is a collar of $\partial\Sigma$ such that 
$\phi|_F=id$. 
Denote $N:=F_\phi\cup (B\times D^2)$, then $V=\pi_tN$ (see (\ref{pit})). 
Let $L_s:=\{ K=s\}\subset N$ denote the $s$-level set of $K$. 
In particular, $L_0=B$. 
Let $i_s:L_s\hookrightarrow N$ denote the inclusion. 
$L_s$ is $S^1_t$-invariant for any $s$. 
Note that $\alpha'|_N$ is independent of $t$, hence 
we have the following simple 
\begin{fact} 
$i_s^*\alpha'$ is a $t$-independent closed 1-form on $L_s$ 
 for  any regular value $s>0$ of $K|_N$. 
\end{fact} 

Let $N_x:=\{ K\leq K|_{t_x}\}\subset N$. 
Let $C\subset M$ be the image of some 
$\tilde{u}\in\tilde{\cM}(\Upsilon,t_x)$.

\begin{lem} 
$C\subset N_x$. 
\end{lem} 

\begin{proof} 
Denote $C_s:=C\cap L_s$ for $L_s\subset N$. 
Suppose that $C\not\subset N_x$. Then there exists a level set 
$L_{s_o}\subset N\setminus N_x$ 
for some $s_o$ such that  $C\pitchfork L_{s_o}$. 
So $C_{s_o}\subset C$ is a union of finitely many embedded circles 
which are pairwise disjoint. Note that 
$[C_{s_o}]=[S^1_t]\in H_1(L_s,\Z)$.  
Let $A\subset C$ denote the domain bounded by $C_{s_o}$ and 
$t_x$ then we have 
\[ 
\int_Ad\alpha '=\int_{t_x}Kdt -
\int_{C_{s_o}}sdt =K(t_x)-s_o<0,  
\] 
which contradicts with (\ref{pos}). So we conclude that 
$C\subset N_x$.  
\end{proof} 

The above lemma implies that, if $\Upsilon\neq\emptyset$ then 
$\Upsilon$ consists of Reeb orbits of $\alpha'$ in 
$N_x\setminus t_x$. Note that $N_x\setminus t_x$ consists of two 
disjoint connected components. We write 
\[ 
N_x\setminus t_x=N_y\cup N_B,  
\] 
where $N_y$ is the connected 
component containing $\bar{t}_y$, and 
$N_B$ is the connected component containing $B$. 

Note that $C$ intersects with $t_x$ at finitely many 
points. If $C\cap t_x\neq\emptyset$ then 
$C\setminus t_x=(C\cap N_y)\cup (C\cap N_B)$ is a union 
of two disjoint set, which is impossible, so $C$ does not 
intersect with $t_x$ geometrically. Therefore we must have 
either $C\subset N_y$ or $C\subset N_B$. 

\begin{lem} \label{CinNB} 
$C\subset N_B$. 
\end{lem} 

\begin{proof} 
Suppose not. Then $C\subset N_y$ and $\Upsilon$ consists of 
Reeb orbits in $N_y$. 
Note that for $K(t_y)<s<K(t_x)$, 
\[ 
L_s\cap N_y\cong S^1_\theta\times S^1_t . 
\] 
Also recall that 
\[ 
\alpha '=\begin{cases} 
h_1(r)dt+h_2(r)d\theta \quad \text{ on } (D_y)_\phi \\ 
Kdt+\beta \hspace{.7in} \text{ on } N_y\setminus (D_y)_\phi 
\end{cases} 
\] 
Then the Reeb vector field $R'$ of $\alpha'$ is 
\[ 
R'=\begin{cases} 
\big( h'_2(r)\partial_t-h'_1(r)\partial_\theta \big) /
(h_1(r)h'_2(r)-h_2(r)h'_1(r)) \quad \text{ on } (D_y)_\phi, \\ 
\big(\partial_t+Y\big) /(K+\beta(Y)) \hspace{1.5in} \text{ on } 
N_y\setminus (D_y)_\phi, 
\end{cases} 
\] 
where $Y$ is the vector field in $\pi_t(N_y)$ such that 
$d\beta(Y,\cdot)=dK$. Parametrize the family of level sets 
$L_s\cap N_y$ by $s$. We have $d\beta =gds\wedge d\theta$ 
for some positive function $g$. Hence
$Y=\frac{-K_s}{g}\partial_\theta$ where $K_s:=\frac{dK}{ds}>0$. 
Therefore if $\gamma\subset N_y$
is a Reeb orbit of $R'$ then  exactly one of the following two
cases must be held: 

\begin{enumerate} 
\item $\gamma$ is a positive iterate of $\bar{t}_y$, hence 
$\gamma$ is homotopic to a curve in $N_y\setminus t_y$ that 
represents the class 
$-n[S^1_t]\in H_1(S^1_\theta\times S^1_t)$ for some 
$n\in \N$; 

\item $\gamma$ is not a positive iterate of $\bar{t}_y$, 
$[\gamma]=-n[S_\theta^1]+m[S_t^1]\in H_1(S^1_\theta\times S^1_t)$ 
for some $n\in N$, $m\in \Z$. 

\end{enumerate} 

On the other hand 
$[t_x]=-[\bar{t}_y]=-[S^1_t]\in H_1(N_y,\Z)$. So 
$\Upsilon \neq\emptyset$ as $t_x$ is not contractible in $N_y$. 
Also $\Upsilon$ does not consists of positive iterates of 
$\bar{t}_y$. Now that $t_x$ is not homologous in 
$N_y\setminus t_y$ to any nonempty finite collection of 
orbits of $R'$, $C$ has to intersect with the Reeb orbit 
$\bar{t}_y$ nontrivially and positively at every point of 
intersection. Let $U\subset N_y$ be a tiny tubular neighborhood 
of $\bar{t}_y$ such that $C\pitchfork \partial U$. 
Let $\sigma:=C\pitchfork \partial U$, then 
$[\sigma]=-n[S^1_\theta]\in H_1(N_y\setminus U)$ for some 
$n\in \N$. Write $C':=C\setminus U$ and 
$\partial C'=\partial_+C'\cup \partial_-C'$, where 
$\partial_+C'=t_x$, $\partial_-C'=\Upsilon\cup \sigma$. 
But $[\Upsilon]=-n'[S^1_\theta]\in H_1(N_y\setminus U,\Z)$ for 
some $n'\in\N$ hence $t_x$ is not homologous in $N_y\setminus U$ 
to $\Upsilon\cup \sigma$, which implies that $C$ does not 
exist if $C\subset N_y$. So $C\subset N_B$. 
\end{proof}

\begin{lem} 
$\tilde{\cM}(\Upsilon,t_x)=\emptyset$ unless $\Upsilon=\emptyset$. 
\end{lem} 

\begin{proof} 
By now we have known that $C\subset N_B$. If $\Upsilon\neq
\emptyset$ then $\Upsilon$ consists of Reeb orbits of $\alpha'$ 
in $N_B$. It is easy to check that $B$ is a generator of 
$H_1(N_B,\Z)\cong \Z$ and every Reeb orbit 
in $N_B$ is homotopic to a positive multiple of $B$, while 
$t_x$ is contractible in $N_B$. So $\Upsilon=\emptyset$. 
\end{proof}

\begin{lem} 
$\tilde{\cM}(t_x)=\cM(t_x)$.  
\end{lem} 

\begin{proof} 
Let $D_x\subset N_B$ be an embedded, $S^1_t$-invariant spanning disc 
of $t_x$ such that $D_x\cap B$ is a point. One can show that 
the {\em Conley-Zehnder index} $CZ(t_x,D_x)$ of $t_x$ relative to 
$D_x$ is 2 (see \cite{RS} and Section 1.2 of \cite{EGH}). 
Note that $t_x$ is not homologically trivial in 
$M\setminus B$. Otherwise there would exist another surface 
$S\subset M\setminus B$ with boundary $\partial S=t_x$. Then 
the closed surface 
$S\cup D_x$ intersects with $B$ at exactly one point, which 
is impossible given the fact that $B$ is homologically trivial. 
Moreover, $H_2(N_B,\Z)=0$ so the Conley-Zehnder index 
$CZ(t_x):=CZ(t_x,D_x)=2$ is independent of the choice of a 
spanning surface of $t_x$. By Lemma \ref{CinNB},  
the fact that $H_2(N_B,\Z)=0$ and the formula for the formal 
dimensions (see \cite{EGH} Proposition 1.7.1 for the general formula)
of components of $\tilde{\cM}(t_x)$ 
we get that $\tilde{\cM}(t_x)$ is of pure dimension with 
\[ 
\dim\tilde{\cM}(t_x)=CZ(t_x)-1=2-1=1 
\]  
provided that $\tilde{\cM}(t_x)\neq\emptyset$. 
Hence $\tilde{\cM}(t_x)=\cM(t_x)$. 

\end{proof}

\begin{lem} 
The algebraic number of $\cM(t_x)/\R$ is $\pm 1$. 
\end{lem}

\begin{proof} 
By Lemma \ref{CinNB}, suppose that $\cM(t_x)\neq\emptyset$ then 
the image $C$ in $M$ of any element $\tilde{u}\subset \cM(t_x)$ 
is contained in $N_B$. Moreover, $C\pitchfork B$ is a single 
point, $C$ intersects with $B$ positively at their point of 
intersection. 
Recall that $S^1_t$ acts on $N_B$ via 
rotations, fixing $B$ pointwise and acting freely on
$N_B\setminus B$. Note that $\alpha'$ is $S^1_t$-invariant. 
To study $\cM(t_x)$ we consider an $\alpha'$-admissible almost 
complex structure $J$ which is also $S^1_t$-invariant. 
Let $\cM^s (t_x)\subset\cM (t_x)$ denote the subset 
consisting of $S^1_t$-invariant elements of $\cM(t_x)$.

Let $i$ denote the standard complex structure on $\C$. 
Let $z=s+it$ be the complex coordinate of $\C$. 
Denote $\tilde{u}_s:=\frac{d\tilde{u}}{ds}$, $\tilde{u}_t
:=\frac{d\tilde{u}}{dt}$. 
A map $\tilde{u}=(a,u):\C\to\R\times M$ is a element of 
$\cM (t_x)$ if $\tilde{u}$ satisfies the d-bar equation 
\begin{equation} \label{d-bar}  
\bar{\partial}_J\tilde{u}:=\tilde{u}_s+J(\tilde{u}) 
\tilde{u}_t=0, 
\end{equation}  
$u(\{ |z|=r\} )\to t_x$ as $r\to \infty$, 
$a(z)\to \infty$ as $z\to \infty$, and $\int_{u(\C)}d\alpha '>0$ 
is finite.  By applying a reparametrization if necessary 
we may assume that $u(0)\in B$.

Most of the proof essentially follows the arguments 
and methods used in Section 7 of \cite{Y2}. Here we only 
outline the idea. 

{\bf Claim 1}: $\cM ^s(t_x)/\R$ consists of a single 
element.  

Like in \cite{Y2}, 
this is done by showing that there is a 1-1 correspondence 
between $\cM^s(t_x)/\R$ and the trajectories from $x$ to $B$ 
of some gradient-like vector field with respect to $-\bar{K}$. 
And there is one and only one such trajectory.

{\bf Claim 2}: For generic $S^1_t$-invariant $J$, 
the linearized d-bar operator is surjective at any 
element of $\cM ^s(t_x)$. 

In particular, this implies that an $S^1_t$-invariant 
solution to the d-bar equation (\ref{d-bar}) 
is an isolated solution. 

Note that 
\begin{equation} \label{J_B} 
T_{\R\times B}(R\times M)=(\underline{\R}\oplus
\underline{R'})\oplus\underline{\C} ,\quad \xi'|_{\R\times B}=
\underline{\C} 
\end{equation}  
where $\underline{\R}$, $\underline{R'}$ denote the trivial 
real line bundles generated by $\partial_t$ and $R'$ the 
Reeb vector field of $\alpha'$ respectively, and $\underline{\C}$ 
is the trivial complex line bundle which equals $\xi'$ when 
restricted to $\R\times B$.

Now that the $\alpha'$-admissible 
almost complex structure $J$ is also $S^1_t$ invariant, then 
\[ 
J|_{\R\times B}=\begin{bmatrix}i&0\\0&i\end{bmatrix}, \quad 
i=\begin{bmatrix}0&-1\\1&0\end{bmatrix},  
\] 
according to the decomposition in (\ref{J_B}). 
 Let $D_{\tilde{u}}$ 
denote the linearized d-bar operator $\bar{\partial}$ at 
$\tilde{u}$. Then for $\eta\in W^{1,2}(\C,\tilde{u}^*T(\R\times
M))$, 
\begin{equation}  \label{linear} 
D_{\tilde{u}}\eta=\eta_s+J\eta_t+\nabla_\eta J.  
\end{equation}  
In particular, near $z=0$ (\ref{linear}) is just the perturbation 
of the standard Cauchy-Riemann equation by a bounded zero order 
term, hence is surjective on $|z|\leq \epsilon$ for some 
$\epsilon >0$.

Now away from $z=0$ we have 
$u:\C\setminus\{ 0\}\cong \R\times S^1\to N_B\setminus B$. 
$N_B\setminus B$ has the structure of a trivial $S^1$ bundle over 
an annulus. Then following Lemma 7.5 of \cite{Y2}, 
for generic $S^1_t$-invariant $J$, $D_{\tilde{u}}$ 
is surjective when $|z|\geq \epsilon$. We then conclude that 
for generic $S^1_t$-invariant $J$, 
$D_{\tilde{u}}$ is surjective at every $S^1_t$-invariant 
solution $\tilde{u}$.

Finally we will show that $\cM^s(t_x)=\cM (t_x)$ essentially. 
This is done by considering branched covers over $N_B$. 

For each $n\in N$ let $\Z_n\subset S^1$ denote the cyclic 
subgroup of order $n$ generated by the $\frac{2\pi}{n}$-rotation 
on $S^1_t$. The action of $\Z_n$ induces an $n:1$ branched 
covering map 
\[ 
\Phi_n:N_B\to N_B \quad \text{ with $B$ the branch set}.
\]  
Each $S^1_t$-invariant $J$ induces an infinite sequence of 
$S^1_t$-invariant almost complex structures 
\[ 
J_n:=(\Phi)_*J(\Phi_n^{-1})_* 
\] 
Let $\alpha'_n:=(\Phi_n)_*\alpha'$. 
$J_n$ is $\alpha'_n$-admissible. Note that 
$\alpha'_n$ is homotopic to 
$\alpha'$ as contact 1-forms, hence 
$\xi'_n:=\ker \alpha'_n$ is isotopic to $\xi'$ 
as contact structures on $N_B$. Also $t_x$ is 
a Reeb orbit of $\alpha'_n$. 

{\bf Claim 3}: For $n$ large enough, $\cM^s_{J_n}(t_x)=
\cM_{J_n}(t_x)$. 

If not, then there is an infinite sequence subsequence $n_i$, 
$n_i\to \infty$ as $i\to\infty$ such that 
$\cM_{J_{n_i}}(t_x)\setminus \cM^s_{J_{n_i}}(t_x)\neq\emptyset$. 
Given $\tilde{v}_i\in  \cM_{J_{n_i}}(t_x)\setminus
\cM^s_{J_{n_i}}(t_x)$, $\tilde{v}_i$ lifts via $\Phi_{n_i}$ 
to an $\Z_{n_i}$-invariant element $\tilde{u}_i\in 
\cM_{J}(t_x)\setminus \cM^s_{J}(t_x)$, $\tilde{u}_i$ is 
unique up to $\R$-translation. By applying $\R$-translations 
if necessary we get that a subsequence of $\tilde{u}_i$, 
also denoted by $\tilde{u}_i$ by the abuse of language, 
will converge to an $S^1_t$-invariant solution $\tilde{u}
\in\cM^s_J(t_x)$ as $i\to\infty$, this contradicts with 
Claim 2, hence is impossible. 

By trading $J$ for $J_n$ for any $n$ large enough we 
have $\cM^s(t_x)=\cM(t_x)$. This complete the proof. 
\end{proof} 

Following the definition of the boundary operator of the contact 
homology we have the following 
 
\begin{lem} 
$\partial t_x=\pm 1$. Hence $H\Theta(M,\xi')=0$. 
\end{lem} 

Thus finished the proof of Theorem \ref{ot}.

\vspace{,2in} 
\noindent 
{\Large \bf Appendix (by Yakov Eliashberg): Sketch of an 
alternative proof} 

\vspace{.2in} 

The following paragraph is a brief description of the argument 
by Eliashberg  of why one can get in the
overtwisted case a contact form with an exactly 1 holomorphic plane
bounded by one of the orbits.

Consider a model of an overtwisted contact structure  which contains a
solid torus with a Lutz $2\pi n$-twist. The contact form in this solid
torus can be chosen as 
$$\alpha=\cos r dz+\sin (  nr)d\varphi $$ where
$0\leq\rho\leq 1$, $z\in\R/\Z$. This is not exactly good formula because
it is not smooth for $r=0$ but can be smoothed without any problems. 
The torus $\rho=\pi/2n$ foliated by horizontal Reeb orbits which bound
holomorphic planes. By making $n$ large we can make the action of these
orbits arbitrarily small. Now it is easy to see explicitly that there is
no other holomorphic planes bounded by these orbits inside the the
considered solid torus. Also there are no other holomorphic curves
different from the planes for which these orbits can be at the positive
end because their action is less than anybody else's action. On the other
hand, if there is a plane bounded by these orbits which goes outside than
the integral of $d\alpha$ along the piece of this curve inside the
solid torus can be made bigger than the action of the orbit. 
This is, of course, Morse-Bott type form, but by a small perturbation we
get 2 orbits out of the whole torus, and one of them has the required
properties.

\end{document}